\documentclass[fleqn]{mat01}
\usepackage{times,mathtimy,amssymb,latexsym}
\begin{document}

\setcounter{page}{375}
\firstpage{375}

\newtheorem{theor}{\bf Theorem}
\newtheorem{lem}{Lemma}

\font\xxx=tibi at 13.5pt
\def\n{\mbox{\xxx{n}}}

\title{The solutions of the $\n$-dimensional Bessel diamond operator and
the Fourier--Bessel transform of their convolution}

\markboth{H\"{u}seyin Y\i ld\i r\i m, M~Zeki Sar\i kaya and Sermin
\"{O}zt\"{u}rk}{Bessel diamond operator}

\author{H\"{U}SEYIN YILDIRIM$^{*}$, M~ZEKI SARIKAYA and SERMIN
\"{O}ZT\"{U}RK}

\address{Department of Mathematics, Faculty of Science and Arts,
Kocatepe University, Afyon, Turkey\\
\noindent $^{*}$Corresponding Author.\\
\noindent E-mail: hyildir@aku.edu.tr; sarikaya@aku.edu.tr; ssahin@aku.edu.tr}

\volume{114}

\mon{November}

\parts{4}

\Date{MS received 21 April 2004; revised 8 September 2004}

\begin{abstract}
In this article, the operator $\Diamond_{B}^{k}$ is introduced and
named as the Bessel diamond operator iterated $k$ times and is defined
by 
\begin{equation*}
\Diamond_{B}^{k} = [ ( B_{x_{1}} + B_{x_{2}} + \cdots +
B_{x_{p}})^{2} - ( B_{x_{p + 1}} + \cdots + B_{x_{p + q}})^{2} ]^{k},
\end{equation*}
where $\;p + q = n, \;B_{x_{i}} = \frac{\partial^{2}}{\partial
x_{i}^{2}} + \frac{2v_{i}}{x_{i}} \frac{\partial}{\partial x_{i}},\;$
where $2v_{i} = 2\alpha_{i} + 1$, $\;\alpha_{i} > - \frac{1}{2}\;$ [8],
$x_{i} > 0$, $i = 1, 2, \ldots, n,\;k$ is a non-negative integer and $n$
is the dimension of $\mathbb{R}_{n}^{+}$. In this work we study the
elementary solution of the Bessel diamond operator and the elementary
solution of the operator $\Diamond _{B}^{k}$ is called the Bessel
diamond kernel of Riesz. Then, we study the Fourier--Bessel transform of
the elementary solution and also the Fourier--Bessel transform of their
convolution.
\end{abstract}


\maketitle

\section{Introduction}

Gelfand and Shilov \cite{2} have first introduced the elementary
solution of the $n$-dimensional classical diamond operator. Later,
Kananthai \cite{3,4,5} has proved the distribution related to the
$n$-dimensional ultra-hyperbolic equation, the solutions of
$n$-dimensional classical diamond operator and Fourier transformation of
the diamond kernel of Marcel Riesz. Furthermore, Kananthai \cite{4} has
showed that the solution of the convolution form $u(x) = (-1)^{k}
S_{2k}(x) \ast R_{2k}(x)$ is a unique elementary solution of the
$\Diamond^{k} u(x) = \delta$.

In this article, we will define the Bessel ultra-hyperbolic type
operator iterated $k$ times with $x\in \mathbb{R}_{n}^{+} = \left\{
x\hbox{:}\ x = (x_{1},\ldots, x_{n}), x_{1}>0, \ldots, x_{n} >
0\right\}$,
\begin{equation*}
\Box_{B}^{k} = (B_{x_{1}} + B_{x_{2}} + \cdots + B_{x_{p}} - B_{x_{p
+ 1}} - \cdots - B_{x_{p + q}})^{k},\quad p + q = n.
\end{equation*}
We will show that the generalized function $R_{2k}(x)$ as defined by
(10) is the unique elementary solution of the operator $\Box
_{B}^{k}$, that is $ \Box _{B}^{k}R_{2k}(x)=\delta $ where $x\in
\mathbb{R}_{n}^{+}$ and $ \delta $ is the Dirac-delta distribution. $S$
is the Shwartz space of any testing functions and $S^{\prime }$ is a space
of tempered distribution.

Furthermore, we will show that the function $E(x)$ as defined by (8) is
an elementary solution of the Laplace--Bessel operator\pagebreak
\begin{equation}
\Delta_{B} = \sum\limits_{i=1}^{n} B_{x_{i}} = \sum\limits_{i=1}^{n}
\left( \dfrac{\partial ^{2}}{\partial x_{i}^{2}} + \dfrac{2v_{i}}{x_{i}} 
\dfrac{\partial }{\partial x_{i}}\right),
\end{equation}
that is, $\Delta_{B}E(x) = \delta$ where $x\in \mathbb{R}_{n}^{+}$.

The operator $\Diamond_{B}^{k}$ can be expressed as the product of the
operators $\Box_{B}$ and $\Delta_{B}$, that is,
\begin{align}
\Diamond_{B}^{k} &= \left[ \left( \sum\limits_{i=1}^{p} B_{x_{i}}
\right)^{2}-\left( \sum\limits_{i=p+1}^{p+q} B_{x_{i}} \right)^{2}
\right]^{k}\nonumber\\[.3pc]
&= \left[ \sum\limits_{i=1}^{p} B_{x_{i}} - \sum\limits_{i=p+1}^{p+q}
B_{x_{i}}\right] ^{k}\left[ \sum\limits_{i=1}^{p}B_{x_{i}} +
\sum\limits_{i= p+1}^{p+q} B_{x_{i}} \right]^{k}\nonumber\\[.3pc]
&= \Box_{B}^{k} \Delta_{B}^{k}.
\end{align}

Denoted by $T^{y}$ the generalized shift operator acting according to
the law \cite{8}
\begin{align*}
\hskip -2cm T_{x}^{y} \varphi (x) &= C_{v}^{\ast } \int_{0}^{\pi}
\!\!\ldots\!\int_{0}^{\pi} \varphi \left( \sqrt{x_{1}^{2} \!+\! y_{1}^{2} \!-\!
2x_{1}y_{1} \cos \theta_{1}}, \ldots, \sqrt{x_{n}^{2} \!+\! y_{n}^{2} \!-\!
2x_{n} y_{n} \cos \theta_{n}} \right)\\[.3pc]
&\quad\, \times \left( \prod\limits_{i=1}^{n} \sin^{2v_{i} - 1} 
\theta_{i} \right) {\rm d}\theta_{1} \ldots {\rm d}\theta_{n},
\end{align*}
where $x, y \in \mathbb{R}_{n}^{+}$, $C_{v}^{\ast} =
\prod_{i=1}^{n} \frac{\Gamma (v_{i}+1)}{\Gamma (\frac{1}{2})
\Gamma (v_{i})}$. We remark that this shift operator is closely
connected with the Bessel differential operator \cite{8}.
\begin{align*}
\frac{{\rm d}^{2}U}{{\rm d}x^{2}} + \frac{2v}{x} \frac{{\rm d}U}{{\rm
d}x} &= \frac{{\rm d}^{2}U}{{\rm d}y^{2}} + \frac{2v}{y} \frac{{\rm
d}U}{{\rm d}y},\\[.3pc]
U(x,0) &= f(x),\\[.3pc]
U_{y}(x,0) &= 0.
\end{align*}
The convolution operator determined by $T^{y}$ is as follows:
\begin{equation}
(f\ast \varphi)(y) = \int_{\mathbb{R}_{n}^{+}} f(y) T_{x}^{y}
\varphi (x) \left( \prod\limits_{i=1}^{n} y_{i}^{2v_{i}} \right) {\rm
d}y.
\end{equation}
Convolution (3) is known as a $B$-convolution. We note the following
properties for the $B$-convolution and the generalized shift operator:

\begin{enumerate}
\renewcommand\labelenumi{(\alph{enumi})}
\item $T_{x}^{y} \cdot 1 =1$.

\item $T_{x}^{0} \cdot f(x)=f(x)$.

\item If $f(x),g(x)$ $\in C(\mathbb{R}_{n}^{+})$, $g(x)$ is a bounded
function, $x>0$ and 
\begin{equation*}
\hskip -1.25pc \int_{0}^{\infty} \left| f\left( x\right) \right| \left(
\prod\limits_{i=1}^{n}x_{i}^{2v_{i}} \right) {\rm d}x < \infty,
\end{equation*}\pagebreak

\noindent then
\begin{equation*}
\hskip -1.25pc \int_{\mathbb{R}_{n}^{+}}\ T_{x}^{y} f(x) g(y) \left(
\prod\limits_{i=1}^{n} y_{i}^{2v_{i}} \right) {\rm d}y =
\int_{\mathbb{R}_{n}^{+}} f(y)\;T_{x}^{y}g(x) \left(
\prod\limits_{i=1}^{n} y_{i}^{2v_{i}} \right){\rm d}y.
\end{equation*}

\item From (c), we have the following equality for $g(x)=1$.
\begin{equation*}
\hskip -1.25pc \int_{\mathbb{R}_{n}^{+}}\ T_{x}^{y} f(x)\left(
\prod\limits_{i=1}^{n} y_{i}^{2v_{i}} \right) {\rm d}y =
\int_{\mathbb{R}_{n}^{+}}\ f(y) \left( \prod\limits_{i=1}^{n}
y_{i}^{2v_{i}} \right) {\rm d}y.
\end{equation*}

\item $(f\ast g)(x)=(g\ast f)(x)$.
\end{enumerate}

The Fourier--Bessel transformation and its inverse transformation are
defined as\break follows \cite{9}:
\begin{align*}
(F_{B}f) (x) &= C_{v}\int_{\mathbb{R}
_{n}^{+}}f(y)\left( \prod\limits_{i=1}^{n}j_{v_{i}-\frac{1}{2}}\left(
x_{i},y_{i}\right) y_{i}^{2v_{i}}\right) {\rm d}y,\\[.3pc]
(F_{B}^{-1}f) (x) &= \left( F_{B}f\right) \left(
-x\right),\;C_{v}=\left( \prod\limits_{i=1}^{n}2^{v_{i}-\frac{1}{2}}\left(
v_{i} + \frac{1}{2}\right) \right)^{-1},
\end{align*}
where $j_{v_{i} - (1/2)} \left( x_{i}, y_{i}\right)$ is the
normalized Bessel function which is the eigenfunction of the Bessel
differential operator. The following are the equalities for 
Fourier--Bessel transformation \cite{6,7}:
\begin{align}
&F_{B}\delta \left( x\right) = 1\nonumber\\[.3pc]
&F_{B}(f\ast g)(x) = F_{B}f(x) \cdot F_{B} g(x).
\end{align}

Now we are finding the solution of the equation 
\begin{equation*}
\Diamond_{B}^{k} u(x) = \sum\limits_{r=0}^{m} \Diamond_{B}^{r}
\delta,\quad \Diamond_{B}^{0}\delta = \delta
\end{equation*}
or 
\begin{equation}
\Box_{B}^{k} \Delta_{B}^{k} u(x) = \sum\limits_{r=0}^{m}
\Box_{B}^{r} \Delta_{B}^{r} \delta.
\end{equation}
In finding the solutions of (5), we use the properties of $B$-convolutions
for the generalized functions.

\begin{lem}
There is the following equality for Fourier--Bessel transformation 
\begin{equation*}
F_{B} ( |x|^{-\alpha}) =2^{n+2\left| v\right|
-2\alpha }\Gamma \left( \frac{n+2\left| v\right| -\alpha }{2}\right) \left[
\Gamma \left( \frac{\alpha }{2}\right) \right] ^{-1}\left| x\right| ^{\alpha
-n-2\left| v\right| },
\end{equation*}
where $\left| v\right| = v_{1} + \cdots + v_{n}$.
\end{lem}

The proof of this lemma is given in \cite{9}.

\begin{lem}
Given the equation $\Delta_{B}E(x) = \delta$ for $x\in
\mathbb{R}_{n}^{+}${\rm ,} where $\Delta_{B}$ is the Laplace--Bessel
operator defined by $(1)${\rm ,}
\begin{equation}
E(x) = -S_{2}(x)\label{5}
\end{equation}
is an elementary solution of the operator $\Delta_{B}$ where 
\begin{equation*}
S_{2}(x) = \frac{2^{n+2\left| v\right| -4}\Gamma \left( \frac{n+2\left|
v\right| -2}{2}\right) }{\prod_{i=1}^{n}2^{v_{i} - \frac{1}{2}}
\Gamma \left( v_{i} + \frac{1}{2}\right)} \left| x\right|^{2-n-2\left|
v\right|}.
\end{equation*}
\end{lem}

\begin{proof}
For the equation $\Delta _{B}E(x)=\delta$, we have 
\begin{equation}
F_{B}\Delta _{B}E=F_{B}\delta.  \label{6}
\end{equation}
First, we consider the left side of (7) 
\begin{align*}
F_{B}\Delta_{B}E(x) &= C_{v}\int_{\mathbb{R}_{n}^{+}}\left(
\Delta _{B}E(y)\right) \left( \prod\limits_{i=1}^{n}J_{v_{i}-\frac{1}{2}
}\left( x_{i}, y_{i}\right) y_{i}^{2v_{i}}\right) {\rm d}y \\[.3pc] 
&= C_{v}\int_{\mathbb{R}_{n}^{+}}\left( \sum\limits_{i=1}^{n}
\frac{\partial ^{2}E(y)}{\partial y_{i}^{2}}\right) \left(
\prod\limits_{i=1}^{n}J_{v_{i}-\frac{1}{2}}\left( x_{i},y_{i}\right)
y_{i}^{2v_{i}}\right) {\rm d}y \\[.3pc]
&\quad\, + C_{v} \int_{\mathbb{R}_{n}^{+}} \left(
\sum\limits_{i=1}^{n} \frac{2v_{i}}{y_{i}} \frac{\partial E(y)}{\partial
y_{i}}\right) \left( \prod\limits_{i=1}^{n} J_{v_{i} - \frac{1}{2}}
\left( x_{i}, y_{i}\right) y_{i}^{2v_{i}} \right) {\rm d}y.
\end{align*}
If we apply partial integration twice in the first integral and once in
the second integral, then we have 
\begin{align*}
F_{B}\Delta_{B}E(x) &= C_{v}\int_{\mathbb{R}_{n}^{+}}E(y)\left(
\sum\limits_{i=1}^{n} \dfrac{\partial^{2}}{\partial y_{i}^{2}}
\prod\limits_{i=1}^{n} J_{v_{i} - \frac{1}{2}}\left( x_{i},y_{i}\right)
\right) y_{i}^{2v_{i}} {\rm d}y \\[.3pc]
&\quad\, + C_{v} \int_{\mathbb{R}_{n}^{+}}\;E(y)\left(
\sum\limits_{i=1}^{n} \frac{2v_{i}}{y_{i}} \frac{\partial }{\partial
y_{i}} \prod\limits_{i=1}^{n} J_{v_{i} - \frac{1}{2}}\left(
x_{i}, y_{i}\right) \right) y_{i}^{2v_{i}} {\rm d}y \\[.3pc]
&= C_{v} \int_{\mathbb{R}_{n}^{+}} E(y) \left(
\sum\limits_{i=1}^{n} B_{y_{i}} \prod\limits_{i=1}^{n} J_{v_{i} -
\frac{1}{2}}\left( x_{i}, y_{i}\right) \right) y_{i}^{2v_{i}} {\rm d}y.
\end{align*}
Here, if we use the following equality \cite{8},
\begin{equation*}
\int_{0}^{\infty} E(y) B_{y_{i}} J_{v_{i} - \frac{1}{2}} \left(
x_{i},y_{i}\right) y_{i}^{2v_{i}} {\rm d} y_{i} = - x_{i}^{2}
\int_{0}^{\infty} E(y) J_{v_{i} - \frac{1}{2}}\left( x_{i},
y_{i}\right) y_{i}^{2v_{i}} {\rm d}y_{i},
\end{equation*}

$\left.\right.$\vspace{-1.5pc}

\noindent then we have
\begin{align*}
F_{B}\Delta_{B}E(x) &= -(x_{1}^{2} + x_{2}^{2} + \cdots + x_{n}^{2})
C_{v}\\
&\quad\, \times \int_{\mathbb{R}_{n}^{+}}E(y)\left( \prod\limits_{i=1}^{n}
J_{v_{i} - \frac{1}{2}} \left( x_{i},y_{i}\right) y_{i}^{2v_{i}} \right)
{\rm d}y \\[.3pc]
&= -\left| x\right| ^{2}F_{B}E(x).
\end{align*}
The right side of equality (7) is $F_{B}\delta =1$. Then 
\begin{equation*}
F_{B}\Delta _{B}E=-\left| x\right| ^{2}F_{B}E=1.
\end{equation*}
From Lemma~1 and inverse Fourier--Bessel transformation we obtain
\begin{equation}
E(x) = -\frac{2^{n+2\left| v\right| -4}\Gamma \left( \frac{n+2\left|
v\right| -2}{2}\right)}{\prod_{i=1}^{n} 2^{v_{i} - \frac{1}{2}}
\Gamma \left( v_{i} + \frac{1}{2}\right) }\left| x\right|^{2-n-2\left|
v\right|}. \label{7}
\end{equation}
That completes the proof.\hfill$\Box$
\end{proof}

\begin{lem}
Given the equation $\Delta _{B}^{k}u(x)=\delta $ for $x\in
\mathbb{R}_{n}^{+}${\rm ,} where $\Delta _{B}^{k}$ is the Laplace--Bessel
operator iterated $k$ times defined by 
\begin{equation*}
\Delta_{B}^{k} = (B_{x_{1}} + B_{x_{2}} + \cdots + B_{x_{n}})^{k},
\end{equation*}
then $u(x)=(-1)^{k}S_{2k}(x)$ is an elementary solution of the operator 
$\Delta _{B}^{k}$ where
\begin{equation}
S_{2k}(x)=\frac{2^{n+2\left| v\right| -4k}\Gamma \left( \frac{n+2\left|
v\right| -2k}{2}\right) }{\prod_{i=1}^{n}2^{v_{i}-\frac{1}{2}}\Gamma
\left( v_{i}+\frac{1}{2}\right) }\left| x\right| ^{2k-n-2\left| v\right|}.
\label{8}
\end{equation}
\end{lem}

The proof of Lemma~3 is similar to the proof of Lemma~2.

\begin{lem}
If $\ \Box_{B}^{k} u(x) = \delta$ for $x\in \Gamma_{+}$ $=\{
x\in \mathbb{R}_{n}^{+}{\rm :}\ x_{1} > 0, x_{2} > 0, \ldots, x_{n} >
0$ and $V > 0\}${\rm ,} where $\Box_{B}^{k}$ is the Bessel-ultra
hyperbolic operator iterated $k$ times defined by
\begin{equation*}
\Box_{B}^{k} = (B_{x_{1}} + B_{x_{2}} + \cdots + B_{x_{p}} -
B_{x_{p+1}} - \cdots - B_{x_{p+q}})^{k},\quad p + q = n,
\end{equation*}

$\left.\right.$\vspace{-1.5pc}

\noindent then $u(x) = R_{2k}(x)$ is the unique elementary solution of the operator 
$\Box_{B}^{k}$ where 
\begin{align}
R_{2k}(x) &= \frac{V^{\left(\frac{2k - n - 2\left| v\right|}{2}\right)}}{K_{n}(2k)}\nonumber\\[.3pc]
&= \frac{(x_{1}^{2} + x_{2}^{2} + \cdots + x_{p}^{2} - x_{p+1}^{2} -
\cdots - x_{p+q}^{2})^{\left(\frac{2k-n-2 \left| v\right| }{2}\right)}}{K_{n}(2k)}
\label{9}
\end{align}
for 
\begin{equation*}
K_{n}(2k) = \frac{\pi^{\frac{n + 2\left| v\right| -1}{2}}\Gamma \left( \frac{
2+2k-n-2\left| v\right| }{2}\right) \Gamma \left( \frac{1-2k}{2}\right)
\Gamma \left( 2k\right) }{\Gamma \left( \frac{2+2k-p-2\left| v\right| }{2}
\right) \Gamma \left( \frac{p-2k}{2}\right)}.
\end{equation*}
\end{lem}

The proof of this lemma can be from Lemmas~1--3.

\begin{lem}
$R_{2k}(x)$ and $S_{2k}(x)$ are homogeneous distributions of order
$(2k-n-2\left| v\right|).$
\end{lem}

\begin{proof}
We need to show that $R_{2k}(x)$ satisfies the Euler equation 
\begin{equation*}
(2k-n-2\left| v\right|) R_{2k}(x) = \sum\limits_{i=1}^{n} x_{i}
\frac{\partial R_{2k}(x)}{\partial x_{i}}.
\end{equation*}

Now 
\begin{align*}
&\sum\limits_{i=1}^{n}x_{i}\dfrac{\partial R_{2k}(x)}{\partial x_{i}}\\[.3pc]
&\quad = \frac{1}{K_{n}(2k)} \sum\limits_{i=1}^{n} x_{i} \frac{\partial }{\partial
x_{i}} (x_{1}^{2} + x_{2}^{2} + \cdots + x_{p}^{2} \!-\! x_{p+1}^{2} \!-
\cdots - x_{p+q}^{2})^{\left(\frac{2k-n-2\left| v\right| }{2}\right)}\\[.3pc] 
&\quad = \frac{(2k-\!n-\!2\left| v\right| )}{K_{n}(2k)} (x_{1}^{2} + x_{2}^{2} +
\cdots + x_{p}^{2} \!-\! x_{p+1}^{2} \!- \cdots -
x_{p+q}^{2})^{\left(\frac{2k-n-2\left| v\right| }{2} - 1 \right)}\\[.3pc] 
&\qquad\, \times (x_{1}^{2} + x_{2}^{2} + \cdots + x_{p}^{2} -
x_{p+1}^{2} - \cdots - x_{p+q}^{2})\\[.3pc]
&\quad = \frac{(2k-n-2\left| v\right| )}{K_{n}(2k)}
(x_{1}^{2} + x_{2}^{2} + \cdots + x_{p}^{2} \!-\! x_{p+1}^{2} \!- \cdots -
x_{p+q}^{2})^{\left(\frac{2k-n-2\left| v\right| }{2}\right)}\\[.3pc] 
&\quad = (2k-n-2\left| v\right| )R_{2k}(x).
\end{align*}
Hence $R_{2k}(x)$ is a homogeneous distribution of order $(2k-n-2\left|
v\right|)$ as required and similarly $S_{2k}(x)$ is also a homogeneous
distribution of order $(2k-n-2\left| v\right|)$.\hfill$\Box$
\end{proof}

\begin{lem}
$R_{2k}(x)$ and $S_{2k}(x)$ are the tempered distributions.
\end{lem}

\begin{proof}
Choose ${\rm supp}\ R_{2k} = K\subset \bar{\Gamma}_{+}$, where $K$ is
a compact set. Then $R_{2k}$ is a tempered distribution with compact
support and by \cite{1}, pp.~156--159, $S_{2k}(x)\ast R_{2k}(x)$ exists
and is a tempered distribution.\hfill$\Box$
\end{proof}

\begin{lem}{\rm (}The $B$-convolutions of tempered distributions{\rm )}
\begin{equation}
S_{2k}(x)\ast S_{2m}(x)=S_{2k+2m}(x).\label{10}
\end{equation}
\end{lem}

\begin{proof}
From Lemmas~1 and 2 we have
\begin{equation*}
F_{B} S_{2k}(x) = - \left| x\right|^{-2k}\quad \hbox{ and }\quad F_{B}
S_{2m} (x) = -\left| x\right| ^{-2m}.
\end{equation*}
Thus from (4) we obtain 
\begin{align*}
F_{B}\left( S_{2k}(x)\ast S_{2m}(x)\right) &= F_{B}S_{2k}(x)F_{B}S_{2m}(x)\\[.3pc]
&= \left| x\right| ^{-2k-2m} \\[.3pc] 
S_{2k}(x)\ast S_{2m}(x) &= F_{B}^{-1}\left| x\right| ^{-2k-2m} \\[.3pc] 
&= C(v,m,k,n)\left| x\right| ^{-2k-2m-n-2\left| v\right| } \\[.3pc] 
&= S_{2k+2m}(x),
\end{align*}
where 
\begin{equation*}
C(v,m,k,n) = \frac{2^{n + 2\left| v\right| - 4(m+k)} \Gamma
\left( \frac{n+2\left| v\right| -2(m+k)}{2} \right)}{\prod_{i=1}^{n}
2^{v_{i} - \frac{1}{2}} \Gamma \left( v_{i} + \frac{1}{2} \right)}.
\end{equation*}

Now from (6) and (11) with $m=k=1$, we have 
\begin{align*}
E(x)\ast E(x) &= (-S_{2}(x))\ast (-S_{2}(x)) \\[.3pc]
&= (-1)^{2}S_{2+2}(x) \\[.3pc] 
&= S_{4}(x).
\end{align*}
By induction, we obtain 
\begin{equation}
\underset{k\ \text{times}}{\underbrace{E(x)\ast E(x)\ast \cdots \ast
E(x)}} = (-1)^{k}S_{2k}(x).\label{11}
\end{equation}
\hfill$\Box$
\end{proof}

\begin{lem}
Given the equation $\Delta _{B}^{k}u(x)=\delta${\rm ,} then 
$u(x)=(-1)^{k}S_{2k}(x)$ is an elementary solution of the operator 
$\Delta_{B}^{k}$ where{\rm ,} $(-1)^{k}S_{2k}(x)$ is defined by $(12)$.
\end{lem}

\begin{proof}
Now $\Delta _{B}^{k}u(x)=\delta$ can be written in the form 
\begin{equation*}
\Delta_{B}^{k}\delta \ast u(x) = \delta.
\end{equation*}
$B$-convolving both sides by the function $E(x)$ defined by (8), we obtain 
\begin{equation*}
(E(x)\ast \Delta _{B}^{k}\delta )\ast u(x)=E(x)\ast \delta =E(x)
\end{equation*}
and 
\begin{equation*}
(\Delta _{B}E(x)\ast \Delta _{B}^{k-1}\delta )\ast u(x)=E(x).
\end{equation*}
Since $\Delta _{B}^{k}E(x)=\delta$ we have 
\begin{equation*}
(\delta \ast \Delta _{B}^{k-1}\delta )\ast u(x)=E(x).
\end{equation*}
Hence 
\begin{equation*}
(\Delta _{B}^{k-1}\delta )\ast u(x)=E(x).
\end{equation*}
By keeping on $B$-convolving $E(x)$, $k-1$ times, we obtain 
\begin{equation*}
\delta \ast u(x)=\underset{k\ \text{times}}{\underbrace{E(x)\ast E(x)\ast
\cdots \ast E(x)}}\ .
\end{equation*}
It follows that 
\begin{equation*}
u(x)=(-1)^{k}S_{2k}(x)
\end{equation*}
by (12) as required.\hfill$\Box$\pagebreak

Before proving the theorems, we need to define the
$B$-convolution of $(-1)^{k}S_{2k}(x)$ with $R_{2k}(x)$ defined by (10)
with $k=0, 1, 2,\ldots\ $. Now for the case $2k\geq n+2\left| v\right|$,
we obtain $(-1)^{k}S_{2k}(x)$ and $R_{2k}(x)$ as analytic functions
that are ordinary functions. Thus the $B$-convolution
\begin{equation}
(-1)^{k}S_{2k}(x)\ast R_{2k}(x)\label{12}
\end{equation}
exists. Now for the case $2k<n+2\left| v\right|$, by Lemma~6 we obtain 
$(-1)^{k}S_{2k}(x)$ and $R_{2k}(x)$ as tempered distributions.

Let $K$ be a compact set and $K\subset \bar{\Gamma}_{+}$ where
$\bar{\Gamma} _{+}$ is defined closer to $\Gamma _{+}.$ Choose the
support of $R_{2k}(x)$ equal to $K,$ then supp$R_{2k}(x)$ is compact
(closed and bounded). So the $B$-convolution
\begin{equation}
(-1)^{k}S_{2k}(x)\ast R_{2k}(x)\label{13}
\end{equation}
exists and is a tempered distribution from Lemma~6.
\end{proof}

\begin{theor}[\!]
Given the equation $\Diamond _{B}^{k}u(x)=\delta$ for $x\in
\mathbb{R}_{n}^{+}${\rm ,} where $\Diamond_{B}^{k}$ is a diamond Bessel
operator iterated $k$ times defined by $(2)$ then
$u(x)=(-1)^{k}S_{2k}(x)\ast R_{2k}(x)${\rm ,} defined by $(13)$ and
$(14)${\rm ,} is a unique elementary solution of the operator 
$\Diamond_{B}^{k}$.
\end{theor}

\begin{proof}
Now $\Diamond _{B}^{k}u(x)=\delta $ can be written as $\Diamond
_{B}^{k}u(x)=\Box _{B}^{k}\Delta _{B}^{k}u(x)$ by (2).
$\Delta_{B}^{k}u(x)=R_{2k}(x)$ is a unique elementary solution of the
operator $\Box_{B}^{k}$ for $n$ odd with $p$ odd and $q$ even,
or for $n$ even with $p$ odd and $q$ odd. By the method of
$B$-convolution, we have 
\begin{equation*}
\Delta _{B}^{k}\delta \ast u(x)=R_{2k}(x).
\end{equation*}
$B$-convolving both sides by $(-1)^{k}S_{2k}(x)$, we obtain 
\begin{equation*}
((-1)^{k} S_{2k}(x)\ast \Delta _{B}^{k}\delta) \ast
u(x)=(-1)^{k}S_{2k}(x)\ast R_{2k}(x)
\end{equation*}
or 
\begin{equation*}
\Delta _{B}^{k} ((-1)^{k} S_{2k}(x)) \ast u(x) = (-1)^{k}
S_{2k}(x) \ast R_{2k}(x).
\end{equation*}
It follows that $u(x) = (-1)^{k}S_{2k}(x)\ast R_{2k}(x)$ by Lemma~8.
That completes the proof.\hfill$\Box$
\end{proof}

\begin{theor}[\!]
For $0 < r < k${\rm ,} 
\begin{equation*}
\Diamond _{B}^{r} ((-1)^{k}S_{2k}(x)\ast R_{2k}(x))
=(-1)^{k-r}S_{2k-2r}(x)\ast R_{2k-2r}(x)
\end{equation*}
and for $k\leq m$
\begin{equation*}
\Diamond_{B}^{m} ((-1)^{k}S_{2k}(x)\ast R_{2k}(x)) =
\Diamond_{B}^{m-k}\delta.
\end{equation*}
\end{theor}

\begin{proof}
From Theorem~1, 
\begin{equation*}
\Diamond_{B}^{k} ((-1)^{k} S_{2k}(x)\ast R_{2k}(x)) =\delta.
\end{equation*}
Thus 
\begin{equation*}
\Diamond_{B}^{k-r} \Diamond_{B}^{r} ((-1)^{k}S_{2k}(x)\ast
R_{2k}(x)) =\delta
\end{equation*}
or 
\begin{equation*}
\Diamond _{B}^{k-r}\delta \ast \Diamond _{B}^{r} ((-1)^{k}
S_{2k}(x) \ast R_{2k}(x)) = \delta.
\end{equation*}
$B$-convolving both sides by $(-1)^{k-r}S_{2k-2r}(x)\ast R_{2k-2r}(x)$,
we obtain 
\begin{align*}
&\Diamond _{B}^{k-r} [(-1)^{k-r}S_{2k-2r}(x)\ast R_{2k-2r}(x)]
\ast \Diamond _{B}^{r} [(-1)^{k}S_{2k}(x)\ast R_{2k}(x)] \\[.3pc] 
&\quad\, = [(-1)^{k-r}S_{2k-2r}(x)\ast R_{2k-2r}(x)] \ast \delta,
\end{align*}
or by Theorem~1 
\begin{equation*}
\delta \ast \Diamond _{B}^{r} [(-1)^{k}S_{2k}(x)\ast
R_{2k}(x)] = [(-1)^{k-r}S_{2k-2r}(x)\ast R_{2k-2r}(x)].
\end{equation*}
It follows that for $0 < r < k$,
\begin{equation*}
\Diamond_{B}^{r} [(-1)^{k}S_{2k}(x)\ast R_{2k}(x)]
= (-1)^{k-r}S_{2k-2r}(x)\ast R_{2k-2r}(x)
\end{equation*}
as required. For $k\leq m$,
\begin{align*}
\Diamond_{B}^{m} ((-1)^{k} S_{2k}(x)\ast R_{2k}(x)) &=  
\Diamond_{B}^{m-k}\Diamond_{B}^{k} [(-1)^{k}S_{2k}(x)\ast R_{2k}(x)]\\[.3pc] 
&= \Diamond_{B}^{m-k}\delta
\end{align*}
by Theorem~1. That completes the proofs.\hfill$\Box$
\end{proof}

\begin{theor}[\!]
Given the linear differential equation 
\begin{equation}
\Diamond _{B}^{k}u(x)=\sum\limits_{r=0}^{m}c_{r}\Diamond _{B}^{r}\delta,
\label{14}
\end{equation}
where the operator $\Diamond _{B}^{k}$ is defined by $(2)${\rm }, $n$ is
odd with $p$ odd and $q$ even{\rm ,} or $n$ is even with $p$ odd and
$q$ odd{\rm ,} $c_{r}$ is a constant{\rm ,} $\delta$ is the Dirac-delta
distribution and $\Diamond_{B}^{0}\delta = \delta$.
\end{theor}

Then the solutions of (15) that depend on the relationship between
the values of $k$ and $m$ are as follows:

\begin{enumerate}
\renewcommand\labelenumi{(\arabic{enumi})}
\item If $m<k$ and $m=0$, then (15) has the solution
$u(x)=c_{0}(-1)^{k}S_{2k}(x)\ast R_{2k}(x),$ which is an elementary
solution of the operator $\Diamond _{B}^{k}$ in Theorem~1 and is the
ordinary function for $2k<n+2\left| v\right|$.

\item If $0<$ $m<k$, then the solution of (15) is 
\begin{equation*}
\hskip -1.25pc u(x) = \sum\limits_{r=1}^{m} [(-1)^{k-r} S_{2k-2r}(x) \ast
R_{2k-2r} (x)]
\end{equation*}
which is an ordinary function for $2k-2r\geq n+2\left| v\right|$, and a
tempered distribution for $2k-2r<n+2\left| v\right|$.

\item If $m\geq k$, and suppose $k\leq n+2\left| v\right| \leq M,$ then
(15) has the solution
\begin{equation*}
\hskip -1.25pc u(x)=\sum\limits_{r=k}^{M}c_{r}\Diamond _{B}^{r-k}
\end{equation*}
which is only the singular distribution.
\end{enumerate}

Proof of this theorem can be easily seen from Theorems~1, 2 and \cite{4}.

\begin{lem}\hskip -.2pc {\rm (}The Fourier--Bessel transformation 
$\Diamond_{B}^{k}\delta${\rm )}.\ \ Let $\| x\| = (x_{1}^{2} + x_{2}^{2}
+ \cdots + x_{n}^{2})^{1/2}$ for $x\in \mathbb{R}_{n}^{+}$. Then 
\begin{equation}
| F_{B}\Diamond _{B}^{k}\delta| \leq C_{\upsilon } \|
x\|^{2k}.\label{16}
\end{equation}
That is{\rm ,} $F_{B}\Delta _{B}^{k}\delta $ is bounded and continuous
on the space $S'$ of the tempered distribution. Moreover{\rm
,} by the inverse Fourier--Bessel transformation
\begin{equation*}
\Diamond _{B}^{k}\delta = C_{\upsilon} F_{B}^{-1} 
[(x_{1}^{2} + x_{2}^{2} + \cdots + x_{p}^{2})^{2} - (
x_{p+1}^{2} + \cdots + x_{p+q}^{2})^{2}]^{k}.
\end{equation*}
\end{lem}

\begin{proof}
From the Fourier--Bessel transform we have 
\begin{align*}
F_{B}\Diamond_{B}^{k}\delta (x) &= C_{v}\int_{\mathbb{R}
_{n}^{+}}\Diamond _{B}^{k}\delta (y)\left( \prod\limits_{i=1}^{n}j_{v_{i} -
\frac{1}{2}}\left( x_{i},y_{i}\right) y_{i}^{2v_{i}}\right) {\rm d}y \\[.3pc] 
&= C_{v}\int_{\mathbb{R}_{n}^{+}}\Delta _{B}^{k}\Box
_{B}^{k}\delta (y)\left( \prod\limits_{i=1}^{n}j_{v_{i}-\frac{1}{2}}\left(
x_{i},y_{i}\right) y_{i}^{2v_{i}}\right) {\rm d}y \\[.3pc] 
&= C_{v}\int_{\mathbb{R}_{n}^{+}}\Delta _{B}^{k}g(y)\left(
\prod\limits_{i=1}^{n}j_{v_{i}-\frac{1}{2}}\left( x_{i},y_{i}\right)
y_{i}^{2v_{i}}\right) {\rm d}y,
\end{align*}
where $g(y)=\Box _{B}^{k}\delta (y).$ For $k\in \mathbb{N}$ we have
\cite{10}
\begin{equation*}
F_{B}(\Delta _{B}^{k})f=(-1)^{k}\left| x\right| ^{2k}F_{B}f.
\end{equation*}
So we have 
\begin{align*}
F_{B}\Diamond _{B}^{k}\delta (x) &= C_{v}(-1)^{k}\left| x\right|^{2k}
F_{B}g(x)\\[.3pc] 
&= C_{v}(-1)^{k} (x_{1}^{2} + \cdots + x_{n}^{2})^{k}
F_{B} \Box_{B}^{k}\delta (x).
\end{align*}
The same way we have following equality:
\begin{equation*}
F_{B} \Box_{B}^{k}\delta (x) = C_{v}(-1)^{k} (x_{1}^{2} +
\cdots + x_{p}^{2} - x_{p+1}^{2} - \cdots - x_{p+q}^{2})^{k}
F_{B} \delta(x).
\end{equation*}
Since $F_{B}\delta (x)=1$, we can write 
\begin{align*}
F_{B}\Diamond _{B}^{k}\delta (x) &= C_{v}(-1)^{2k} (x_{1}^{2} +
\cdots + x_{n}^{2})^{k}\\[.3pc]
&\quad\, \times (x_{1}^{2} + \cdots + x_{p}^{2} -
x_{p+1}^{2} - \cdots - x_{p+q}^{2})^{k}\\[.3pc]
&= C_{v} [(x_{1}^{2} + \cdots + x_{p}^{2})^{2} - 
(x_{p+1}^{2} + \cdots + x_{p+q}^{2})^{2}]^{k}.
\end{align*}
Then there is the following inequality:
\begin{align*}
|F_{B} \Diamond_{B}^{k}\delta| &= C_{v} [| x_{1}^{2} + \cdots +
x_{n}^{2}| |x_{1}^{2} + \cdots + x_{p}^{2} - x_{p+1}^{2} - \cdots -
x_{p+q}^{2}|]^{k}\\[.3pc]
&\leq C_{v} [|x_{1}^{2} + \cdots + x_{n}^{2}|]^{2k}\\[.3pc]
&= C_{v} \| x\| ^{4k}.
\end{align*}
Therefore $F_{B}\Diamond_{B}^{k}$ is bounded and continuous on the space 
$S^{\prime }$ of the tempered distribution.\pagebreak

Since $F_{B}$ is 1-1 transformation from $S^{\prime }$ to 
$\mathbb{R}_{n}^{+}$, there is the following equation:
\begin{equation*}
\Diamond_{B}^{k}\delta = C_{\upsilon} F_{B}^{-1} [(x_{1}^{2}
+ x_{2}^{2} + \cdots + x_{p}^{2})^{2} - (x_{p+1}^{2} +
\cdots + x_{p+q}^{2})^{2}]^{k}.
\end{equation*}
That completes the proof.\hfill$\Box$
\end{proof}

\begin{theor}[\!]
\begin{align*}
&F_{B} [(-1)^{k} S_{2k}(x) \ast R_{2k}(x)]\\[.3pc]
&\quad\ \leq \frac{C_{\upsilon }}{[(x_{1}^{2} + x_{2}^{2} + \cdots +
x_{p}^{2})^{2} - (x_{p+1}^{2} + \cdots + x_{p+q}^{2})^{2}]^{k}}
\end{align*}
and 
\begin{equation*}
|F_{B} [(-1)^{k} S_{2k}(x)\ast R_{2k}(x)]|
= C_{\upsilon }M\quad \text{for\ a\ large} \  x_{i}\in \mathbb{R}^{+},
\end{equation*}
where $M$ is a constant. That is{\rm ,} $F_{B}$ is bounded and
continuous on the space $S'$ of the tempered distribution.
\end{theor}

\begin{proof}
By Lemma~8, 
\begin{equation*}
\Diamond_{B}^{k} [(-1)^{k}S_{2k}(x)\ast R_{2k}(x)] =\delta
\end{equation*}
or 
\begin{equation}
(\Diamond_{B}^{k}\delta)\ast [(-1)^{k}S_{2k}(x)\ast
R_{2k}(x)] = \delta. \label{17}
\end{equation}
If we applied the Fourier--Bessel transform on both sides of (17), then
we obtain 
\begin{align*}
&F_{B} [(\Diamond_{B}^{k} \delta)\ast [(-1)^{k} S_{2k}(x)\ast
R_{2k}(x)]] = F_{B}\delta (x)\\[.3pc]
&C_{v} \left\langle (\Diamond _{B}^{k}\delta )\ast [(-1)^{k}
S_{2k}(x) \ast R_{2k}(x)], \ \prod\limits_{i=1}^{n} j_{v_{i} -
\frac{1}{2}} (x_{i}, y_{i}) y_{i}^{2v_{i}} \right\rangle =
C_{v}.
\end{align*}
By the properties of $B$-convolution 
\begin{align*}
&C_{v} \left\langle (\Diamond _{B}^{k}\delta), \left\langle 
[(-1)^{k} S_{2k}(x) \ast R_{2k}(x)], \phantom{\prod\limits_{i=1}^{n}} \right.\right.\\[.3pc]
&\quad\, \left.\left. \prod\limits_{i=1}^{n}
j_{v_{i} - \frac{1}{2}} \left( z_{i}, y_{i} \right) y_{i}^{2v_{i}}
\prod\limits_{i=1}^{n} j_{v_{i} - \frac{1}{2}} \left( x_{i}, y_{i}
\right) y_{i}^{2v_{i}} \right\rangle \right\rangle = C_{v},\\[.3pc]
&C_{v} \left\langle [ (-1)^{k} S_{2k}(x) \ast R_{2k}(x)], \phantom{\prod\limits_{i=1}^{n}} \right.\\[.3pc]
&\quad\,\left. \prod\limits_{i=1}^{n} j_{v_{i} - \frac{1}{2}} \left( z_{i}, y_{i} \right)
y_{i}^{2v_{i}} \right\rangle \left\langle (\Diamond_{B}^{k} \delta),
\prod\limits_{i=1}^{n} j_{v_{i} - \frac{1}{2}} \left( x_{i}, y_{i}
\right) y_{i}^{2v_{i}} \right\rangle = C_{v},\\[.3pc]
&F_{B} [(-1)^{k} S_{2k}(x) \ast R_{2k}(x)] \frac{1}{C_{v}}
F_{B} (\Diamond_{B}^{k} \delta) = C_{v}.
\end{align*}
By Lemma~9,
\begin{align*}
&F_{B} [(\Diamond_{B}^{k} \delta)\ast [(-1)^{k} S_{2k}(x)
\ast R_{2k}(x)]]\\[.3pc]
&\quad\ \times [(x_{1}^{2} + x_{2}^{2} + \cdots + x_{p}^{2})^{2} -
(x_{p+1}^{2} + \cdots + x_{p+q}^{2})^{2}]^{k} = C_{v}.
\end{align*}
It follows that 
\begin{align*}
&F_{B} [(-1)^{k} S_{2k}(x) \ast R_{2k}(x)]\\[.3pc]
&\quad\, = \frac{C_{\upsilon}}{[(x_{1}^{2} + x_{2}^{2} + \cdots +
x_{p}^{2})^{2} - (x_{p+1}^{2} + \cdots + x_{p+q}^{2})^{2}]^{k}}.
\end{align*}
Now 
\begin{align}
&|F_{B} [(-1)^{k} S_{2k}(x) \ast R_{2k}(x)]|\nonumber\\[.3pc]
&\quad\, = \frac{C_{v}}{|x_{1}^{2} + x_{2}^{2} + \cdots + 
x_{n}^{2}|^{k} |x_{1}^{2} + \cdots + x_{p}^{2} - x_{p+1}^{2} - \cdots
- x_{p+q}^{2}|^{k}},\label{18}
\end{align}
where $x = (x_{1}, \ldots, x_{n}) \in \Gamma_{+}$ with $\Gamma_{+}$
defined by Lemma~4. Then $( x_{1}^{2} + \cdots + x_{p}^{2} -
x_{p+1}^{2}$ $- \cdots - x_{p+q}^{2}) > 0$ and for a large $x_{i}$
and a large $k$, the right-hand side of (18) tends to zero. It follows
that it is bounded by a positive constant say $M$, that is, we obtain
(17) as required and also by (17). $F_{B}$ is continuous on the space
$S^{\prime }$ of the tempered distribution.\\
$\left.\right.$\hfill$\Box$
\end{proof}

\begin{theor}[\!]
\begin{align*}
&F_{B} [[(-1)^{k} S_{2k}(x) \ast R_{2k}(x)] \ast 
[(-1)^{k} S_{2m}(x) \ast R_{2m}(x)]] \\[.3pc] 
&\quad\, = F_{B} [(-1)^{k}S_{2k}(x)\ast R_{2k}(x)] F_{B} 
[(-1)^{m}S_{2m}(x)\ast R_{2m}(x)]\\[.3pc]
&\quad\, = \frac{C_{\upsilon }}{(( x_{1}^{2} + x_{2}^{2} + \cdots +
x_{p}^{2})^{2} - (x_{p+1}^{2} + \cdots + x_{p+q}^{2})^{2})^{k+m}},
\end{align*}
where $k$ and $m$ are non-negative integers and $F_{B}$ is bounded and
continuous on the space $S'$ of the tempered distribution.
\end{theor}

\begin{proof}
Since $S_{2k}(x)$ and $R_{2k}(x)$ are tempered distribution with compact
support, from Lemmas~6 and 7 we have
\begin{align*}
&[(-1)^{k}S_{2k}(x) * R_{2k}(x)] * [(-1)^{m} S_{2m}(x) *
R_{2m}(x)]\\[.3pc]
&\quad\, = (-1)^{k+m} [S_{2k}(x) * S_{2m}(x)] * [R_{2k}(x) *
R_{2m}(x)]\\[.3pc]
&\quad\, = (-1)^{k+m} S_{2(k+m)}(x) * R_{2(k+m)}(x).
\end{align*}
Taking the Fourier transform on both sides and using Theorem 4 we obtain
\begin{align*}
&F_{B} [(-1)^{k}S_{2k}(x) * R_{2k}(x)] * [(-1)^{m}S_{2m}(x) *
R_{2m}(x)]\\[.3pc]
&\quad\, =\frac{C_{v}}{[x_{1}^{2} + x_{2}^{2} + \cdots + x_{p}^{2})^{2} -
(x_{p+1}^{2} + \cdots + x_{p+q}^{2})^{2}]^{k+m}}
\end{align*}
\begin{align*}
&= \frac{1}{C_{v}} \frac{C_{v}}{[(x_{1}^{2} + x_{2}^{2} + \cdots +
x_{p}^{2})^{2} - (x_{p+1}^{2} + \cdots + x_{p+q}^{2})^{2}]^{k}}\\[.3pc]
&\qquad\,\times \frac{C_{v}}{[(x_{1}^{2} + x_{2}^{2} + \cdots + x_{p}^{2})^{2} -
(x_{p+1}^{2} + \cdots + x_{p+q}^{2})^{2}]^{m}}\\[.3pc]
& = \frac{1}{C_{v}} F_{B} [(-1)^{k} S_{2k}(x) * R_{2k}(x)] F_{B} [(-
1)^{m}S_{2m}(x) * R_{2m}(x)].
\end{align*}

Since $(-1)^{k+m}S_{2(k+m)}(x) * R_{2(k+m)}(x) \in S'$, the space of
tempered distribution, and Theorem 4 we obtain that $F_{B}$ is bounded
and continuous on $S'$.\hfill$\Box$
\end{proof}

\end{document}